\documentclass[12pt]{amsart}
\usepackage{amssymb}
\usepackage{amsthm}
\usepackage{amscd}

\usepackage{amsmath}
\usepackage{epic,eepic}
\usepackage{graphicx}
\DeclareMathAlphabet{\mathpzc}{OT1}{pzc}{m}{it}

\usepackage{latexsym}

\usepackage{pifont}

\theoremstyle{plain}
\newtheorem{lemma}{Lemma}[section]
\newtheorem{prop}[lemma]{Proposition}
\newtheorem{thm}[lemma]{Theorem}
\newtheorem{cor}[lemma]{Corollary}
\newtheorem{aplemma}{Lemma~A.\hspace{-1.5mm}}
\newtheorem{approp}{Proposition~A.\hspace{-1.5mm}}
\newtheorem{apthm}{Theorem~A.\hspace{-1.5mm}}
\newtheorem{apcor}{Corollary~A.\hspace{-1.5mm}}
\newtheorem{intthm}{Theorem}

\usepackage[all]{xy}

\newtheorem{conj}[lemma]{Conjecture}

\theoremstyle{definition}

\newtheorem{rema}{Remark}

\newtheorem{remb}{Remark}

\newtheorem{defi}[lemma]{Definition}
\newtheorem{exa}[lemma]{Example}
\newtheorem{aprem}{Remark~A.\hspace{-1.5mm}}
\newtheorem{apdefi}{Definition~A.\hspace{-1.5mm}}
\newcommand{\bde}{\begin{defi}}
\newcommand{\ede}{\end{defi}\vspace{1mm}}
\newcommand{\ble}{\begin{lemma}}
\newcommand{\ele}{\end{lemma}}
\newcommand{\bpr}{\begin{prop}}
\newcommand{\epr}{\end{prop}}
\newcommand{\bt}{\begin{thm}}
\newcommand{\et}{\end{thm}}
\newcommand{\bco}{\begin{cor}}
\newcommand{\eco}{\end{cor}}
\newcommand{\bre}{\begin{rema}}
\newcommand{\ere}{\end{rema}}
\newcommand{\brea}{\begin{rema}}
\newcommand{\erea}{\end{rema}\vspace{1mm}}
\newcommand{\breb}{\begin{remb}}
\newcommand{\ereb}{\end{remb}\vspace{1mm}}
\newcommand{\bex}{\begin{exa}}
\newcommand{\eex}{\end{exa}}
\newcommand{\bpf}{\begin{proof}}
\newcommand{\epf}{\end{proof}\vspace{1mm}}

\newcommand{\bade}{\begin{apdefi}}
\newcommand{\eade}{\end{apdefi}}
\newcommand{\bale}{\begin{aplemma}}
\newcommand{\eale}{\end{aplemma}}
\newcommand{\bapr}{\begin{approp}}
\newcommand{\eapr}{\end{approp}}
\newcommand{\bat}{\begin{apthm}}
\newcommand{\eat}{\end{apthm}}
\newcommand{\baco}{\begin{apcor}}
\newcommand{\eaco}{\end{apcor}}
\newcommand{\bare}{\begin{aprem}}
\newcommand{\eare}{\end{aprem}}

\newcommand{\be}{\begin{enumerate}}
\newcommand{\ee}{\end{enumerate}}
\newcommand{\bcd}{\[\begin{CD}}
\newcommand{\ecd}{\end{CD}\]}
\newcommand{\bit}{\begin{itemize}}
\newcommand{\eit}{\end{itemize}}
\newcommand{\bq}{\begin{quote}}
\newcommand{\eq}{\end{quote}}
\newcommand{\ba}{\begin{array}}
\newcommand{\ea}{\end{array}}

\newcommand{\mcB}{\mathcal{B}}

\newcommand{\mcL}{\mathcal{L}}

\newcommand{\mcO}{\mathcal{O}}

\newcommand{\mbF}{\mathbb{F}}

\newcommand{\mbP}{\mathbb{P}}
\newcommand{\mbQ}{\mathbb{Q}}

\newcommand{\migi}{\rightarrow}

\newcommand{\isom}{\stackrel{\sim}{\migi}}

\newcommand{\migiincl}{\hookrightarrow}

\newcommand{\migisurj}{\twoheadrightarrow}

\newcommand{\mr}{\mathrm}
\newcommand{\hidden}[1]{\,}

\newcommand{\I}{\u{\i}}

\begin{document}

\title[An effective version of Bely\I's theorem in positive characteristic]{An effective version of Bely\I's theorem \\ in positive characteristic}
\author{Yasuhiro Wakabayashi}
\date{}
\markboth{Yasuhiro Wakabayashi}{}
\maketitle
\footnotetext{Y. Wakabayashi: Department of Mathematics, Tokyo Institute of Technology, 2-12-1 Ookayama, Meguro-ku, Tokyo 152-8551, JAPAN;}
\footnotetext{e-mail: {\tt wkbysh@math.titech.ac.jp};}
\begin{abstract}
The purpose of the present paper is to give an effective version of the noncritical $p$-tame Bely\u{\i} theorem. That is to say, we compute explicitly an upper bound of the minimal degree of tamely ramified Bely\u{\i} maps in positive characteristic which are unramified at a prescribed finite set of points.

\end{abstract}
\tableofcontents 

\section{Introduction: 
Effective noncritical $p$-tame Bely\I \ theorem}
\vspace{5mm}

Let $X$ be a curve over a field $k$, which means, in the present paper, a geometrically connected proper smooth scheme over  $k$ of dimension $1$.
A {\bf Bely\I \ map} on  $X$  is a dominant $k$-morphism  $f : X \migi \mbP^1_k$ with $\mr{Br}(f) \subseteq \{ 0,1,\infty\}$, where $\mbP^1_k$ denotes the projective line over $k$ and $\mr{Br}(f)$
  denotes the set of branch points of $f$.
A celebrated theorem of Bely\I  \ (cf. ~\cite{Bel1}, ~\cite{Bel2}) asserts that
if $k = \overline{\mbQ}$ (i.e., an algebraic closure of the field of rational numbers), then $X$ always admits at least one Bely\I \ map.
 This result has attracted much attention ever since Grothendieck noticed in his  article (cf. ~\cite{Gro}) that it implies surprising correspondences between curves defined  over number fields and  a certain class of bipartite  graphs embedded   in a topological surface called dessins d'enfants.

There are several variations and enhancements of Bely\I's theorem.
For example, Mochizuki (cf. ~\cite{Moc}) and Scherr-Zieve (cf. ~\cite{SZ}) proved the noncritical enhancement, asserting the existence of Bely\I \ maps which are unramified at a prescribed finite set of points.
Khadjavi (cf. ~\cite{Kha}) and Li\c{t}canu (cf. ~\cite{Lit})
 considered the effective versions, which give explicit upper bounds for the minimal degree  of such morphisms  
 in terms of  height.
 The bounds can be interpreted, for example,
bounding the number of edges of a dessin d'enfant given by
the corresponding Bely\I \ map (cf. ~\cite{Gro}). 
On the other hand, we can find analogues of Bely\I's theorem in positive characteristic, i.e., the {\it $p$-tame Bely\I \ theorem} (asserting the existence of a tamely ramified Bely\I \ map) and the {\it $p$-wild Bely\I \ theorem} (asserting the existence of a Bely\I \ map  which admits at most one branch point).
See 
~\cite{Anb}, ~\cite{Ful}, ~\cite{Ka}, ~\cite{Sai}, ~\cite{SY}, and ~\cite{Zap}.

The purpose of the present paper is to give
an effective version of the noncritical $p$-tame Bely\I \ theorem in positive characteristic.
(The  $p$-wild case can be obtained immediately from the previous works, see Proposition \ref{P044})
A  point of our study is that although the effective bounds of Bely\I \ maps
  obtained so far have been given only for $X = \mbP^1_{\overline{\mbQ}}$,
   we compute, in our situation (i.e., the case of positive characteristic),  an upper bound for an arbitrary curve.
 We shall state the main theorems.
Denote by $\mbF_q$ the finite field with $q$ elements, where $q$ is   a power of an odd prime $p$, and by
$\overline{\mbF}_q$ its algebraic closure.
Let  $X$ be a curve over $\mbF_q$ and
let 
$S$,  $T$ 
be (possibly empty) finite sets of $\mbF_q$-rational points of $X$ with $S \cap T = \emptyset$, where $s := \sharp S$, $t:= \sharp T$.
For each field $k$ over $\mbF_q$, we write
$X_k$, $S_k$, and $T_k$ for  the base-changes over $k$ of
$X$, $S$, and $T$ respectively.
Here, by a {\bf $p$-tame Bely\I \ map} on $(X, S, T)$ over $k$, we shall mean 
a tamely ramified $k$-dominant morphism $f : X_{k} \migi \mbP^1_{k}$  satisfying the following conditions:
\begin{align}
f(S_{k}) \cup \mr{Br}(f) \subseteq \{0,1,\infty \}, \hspace{10mm} \{ 0,1,\infty \} \cap f(T_{k}) = \emptyset.
\end{align}
The {\bf $p$-tame Bely\I \ degree} of $(X, S, T)$ is defined as
\begin{align}
{^t \mcB} (X, S, T) := \mr{min} \left\{ \mr{deg}(\phi) \, | \, \text{$\phi$ is  a $p$-tame Bely\I \ map $\phi$ on $(X,S, T)$ over $\overline{\mbF}_q$}\right\}
\end{align}
(where ${^t \mcB} (X, S, T):= \infty$ if there is no $p$-tame Bely\I \ map on $(X, S, T)$ over $\overline{\mbF}_q$).
The value ${^t \mcB} (X) := {^t \mcB} (X, \emptyset, \emptyset)$ (i.e., the minimal degree of tamely ramified Bely\I \ maps on $X$) is simply referred as the {\bf $p$-tame Bely\I \ degree} of $X$.
Then,  the main assertion of the present paper is as follows.
(Note that the upper bound of the degree asserted in the theorem is somewhat rough; by treating strictly various inequalities at each step in our proof, e.g., the inequalities appearing in Proposition \ref{P02},   we can obtain a  sharper bound.)

\vspace{3mm}
\begin{intthm}[Effective version of the noncritical $p$-tame Bely\I \ theorem]
\leavevmode\\
 \ \ \ 
There exists at least one $p$-tame Bely\I \ map on $(X, S, T)$ over $\overline{\mbF}_q$, and moreover, the following inequality holds:
\begin{align}
{^t \mcB} (X, S, T) \leq (2g+t+1) \cdot (q^{\lceil\mr{log}_q (10^2 \cdot (2g+t+1)! \cdot (2g+t+s+1)^2 \cdot \left(\frac{5}{6}\right)^{2g+t+1}\rceil \cdot L (6g+2t)}-1)^{6g+s+2t+1}.
\end{align}
Here, $\lceil  - \rceil$ denotes the ceiling function and, for each nonnegative integer $m$, $L(m)$ denotes the least common multiple of  $1,2, \cdots, m$ (where $L (0) :=1$ if $m=0$).
In particular, if $S=T= \emptyset$, then
we have
\begin{align}
{^t \mcB} (X) \leq (2g+1) \cdot (q^{\lceil\mr{log}_q (10^2 \cdot (2g+1)! \cdot (2g+1)^2 \cdot \left(\frac{5}{6}\right)^{2g+1}\rceil \cdot L (6g)}-1)^{6g+1}
\end{align}
 \end{intthm}
\vspace{3mm}

\vspace{5mm}
\hspace{-4mm}{\bf Acknowledgement} \leavevmode\\
 \ \ \ 
The author would like to thank 
all algebraic curves in positive characteristic, living in the world of mathematics, for their useful comments and heartfelt encouragement.
The author was partially  supported by the Grant-in-Aid for Scientific Research (KAKENHI No.\,18K13385).

\vspace{10mm}
\section{First step of the proof} \vspace{3mm}

%


In this section, we prove Proposition \ref{P02} described below, which is a slightly strengthened version of ~\cite{Ful}, Proposition 8.1, and also is  the main body of Theorem A. (The statement of Proposition \ref{P02} include the case where the base field is algebraically closed, but this will not be used in our proof of  Theorem A.)
Here, recall that a dominant morphism $f: Y \migi \mbP^1_k$ from a curve $Y$ over a field $k$ onto $\mbP^1_k$
 is called a {\it simple covering} if
 the discriminant $\delta (f)$ of $f$ is a simple divisor, i.e.,  has no multiple components after base-change over  an algebraically closed field $\overline{k}$ over $k$.
(If $f$ is as above, then  $\mr{deg}(\delta (f)) = 2 g + 2  \mr{deg}(f)-2$, where $g$ denotes the genus of $Y$.)
In particular, if  $k$ has characteristic $\neq 2$, then
 any simple covering is tamely ramified (cf. ~\cite{Ful}, Theorem 5.6).

\vspace{3mm}
\bpr \label{P02} \leavevmode\\
 \ \ \
 Let $k$ be a field of characteristic $\neq 2$,  
 $X$  
 a curve over  $k$ of genus $g \ (\geq 0)$, and
 $S$, $T$  (possibly empty) finite sets of $k$-rational points of $X$ with $S \cap T = \emptyset$.
  Write $s = \sharp S$ and $t := \sharp T$.
 Also, let $n$ be a positive integer with $n \geq  g +  \mr{max} \{t, g\}$, and
 suppose that one of the following conditions $(*)$, $(**)$  is satisfied:
\begin{itemize}
\item[$(*)$]
$k$ is algebraically closed; 
\vspace{2mm}
\item[$(**)$]
$k=\mbF_q$ for a power $q$ of a prime $p$ 
such that there exists an integer $A$ with $A \geq 3$ satisfying the  inequality:
\begin{align} \label{Erfhih}
q \geq \mr{max} \left\{A^2 \cdot g^2,  10^2 \cdot n! \cdot (n^2 + s) \cdot \left(\frac{5A+4}{9A-6} \right)^{n} \right\}.
\end{align}
\end{itemize}
 Then, there exists a simple covering $\zeta : X  \migi \mbP^1_{k}$ of degree $n$ such that $\zeta (S) \cap \zeta (T) = \emptyset$ and  $\zeta (T)$ consists of one point (if $T \neq \emptyset$) over which $\zeta$ is unramified.
 In particular,  if $k=\mbF_q$
 with $q \geq 10^2 \cdot (2g+t+1)! \cdot (2g+t+s+1)^2 \cdot (5/6)^{2g+t+1}$, then (by considering the case where $n = 2g+t+1$ and $A=3$, we see that)
 there exists a simple covering $\zeta : X \migi \mbP^1_{\mbF_q}$ of degree $2g+t+1$ satisfying the above requirements.  
 \epr
\begin{proof}

For each positive integer $r$, we shall denote by $X^{(r)}$ the $r$-th symmetric product of $X$ over $k$;
it is a geometrically connected proper smooth scheme over $k$ of dimension $r$, each of whose $k$-rational  point corresponds to
an effective divisor of degree $r$ on $X$.
If, moreover, $X^r$ denotes the product of $r$ copies of $X$ over $k$, then we have the natural projection
\begin{align} \label{E030}
\pi_r : X^r \migi X^{(r)}.
\end{align} 
Let us first consider the following lemma.

\vspace{3mm}
\ble \label{L005} \leavevmode\\
 \ \ \ 
Suppose that $(**)$ is satisfied.
\begin{itemize}
\item[(i)]
There exists a subset $R$ of
 $X (\mbF_p) \setminus (S \cup T)$ with $\sharp R = n-t-g$.
\vspace{1mm}
\item[(ii)]
For each positive integer $r$, the following inequalities hold:
\begin{align}
 \frac{1}{7 \cdot r!} \cdot \left( 3-\frac{2}{A}\right)^r \cdot q^r < \sharp X^{(r)}(\mbF_q) <\left(\frac{5}{3}+ \frac{4}{3A} \right)^{r} \cdot q^r.
\end{align}
\end{itemize}
\ele
\begin{proof}
Recall the Hasse-Weil theorem, asserting the following equality for each positive integer $m$: 
\begin{align} \label{E010}
| \sharp X (\mbF_{q^m}) - (q^m+1)| \leq 2g \sqrt{q^m}.
\end{align}
Then, assertion (i) follows from the above inequality for $m=1$ and the inequality  $q +1 - 2g \sqrt{q} \geq n+s-g \ \left(= (n-t-g)+s+t\right)$ arising immediately from the assumption
(\ref{Erfhih}).

In what follows, let us prove assertion (ii).
The assumption
 $q \geq A^2 \cdot g^2$ ($\Longrightarrow  \frac{2}{A} \cdot q^m \geq 2 g \sqrt{q^m}$ for any $m \geq 1$) and (\ref{E010}) imply
the inequalities:
\begin{align} \label{E012}
\left(1-\frac{2}{A}\right)\cdot q^m+1 \leq \sharp X (\mbF_{q^m})\leq \left(1+ \frac{2}{A}\right) \cdot q^m+1.
\end{align}
Since $X^{(r)}(\mbF_q)$ is the set of effective divisors of degree $r$ on $X$, 
its number can be calculated by
\begin{align} \label{E011}
\sharp X^{(r)} (\mbF_q) = \sum_{i=1}^r \frac{1}{i!}\sum_{\genfrac{.}{.}{0pt}{}{a_1, \cdots, a_i \geq 1}{\sum_{j=1}^i a_j =r}} \prod_{j=1}^i \frac{\sharp X (\mbF_{q^{a_j}})}{a_j}. 
\end{align}
By means of  (\ref{E011}) and the first  inequality of (\ref{E012}),
we have the following sequence of inequalities:
\begin{align}
\sharp X^{(r)}(\mbF_q)
& \geq 
\sum_{i=1}^r \frac{1}{i!}\sum_{\genfrac{.}{.}{0pt}{}{a_1, \cdots, a_i \geq 1}{\sum_{j=1}^i a_j =r}} \prod_{j=1}^i 
\frac{1}{a_j} \cdot \left( \left(1-\frac{2}{A} \right) \cdot q^{a_j} +1\right) 
 \\
& > \sum_{i=1}^r \frac{1}{i!}\sum_{\genfrac{.}{.}{0pt}{}{a_1, \cdots, a_i \geq 1}{\sum_{j=1}^i a_j =r}}  
\left(\frac{r}{\sum_{j=1}^i a_j}\right)^r \prod_{j=1}^i\left(1-\frac{2}{A}\right) \cdot q^{a_j} \notag \\
& = q^r \cdot \sum_{i=1}^r \frac{1}{i!}\sum_{\genfrac{.}{.}{0pt}{}{a_1, \cdots, a_i \geq 1}{\sum_{j=1}^i a_i =r}} \left(1-\frac{2}{A} \right)^i \notag \\
&  =
   q^r \cdot \left(1- \frac{2}{A} \right) \cdot \sum_{i=1}^r \frac{1}{i!} \cdot  \left(1- \frac{2}{A} \right)^{i-1} \cdot \binom{r-1}{i-1} \notag  \\
   & \geq  q^r  \cdot \left(1- \frac{2}{A} \right) \cdot \sum_{i=1}^r \frac{2^{(r-1)-(i-1)}}{r!}  \cdot  \left(1- \frac{2}{A} \right)^{i-1} \cdot \binom{r-1}{i-1} \notag \\
   & \geq 
 q^r \cdot \frac{1}{r!} \cdot \left(1-\frac{2}{A} \right) \cdot \left(3 - \frac{2}{A} \right)^{r-1} \notag \\
 & \geq q^r \cdot \frac{1}{7 \cdot r!} \cdot \left( 3-\frac{2}{A}\right)^r, \notag
\end{align}
where the 
second inequality follows from the geometric-harmonic inequality 
$\left(\prod_{j=1}^i \frac{1}{a_j}\right)^{1/r} \geq \frac{r}{\sum_{j=1}^i a_j}$.
On the other hand, 
let $a := \frac{3}{5}$ if $A=3$ and $a:= \frac{2}{3}$ if $A\geq 4$.
Then,
by (\ref{E011}) and the second inequality of (\ref{E012}),
we have
\begin{align}
\sharp X^{(r)}(\mbF_q)  
& \leq 
\sum_{i=1}^r \frac{1}{i!}\sum_{\genfrac{.}{.}{0pt}{}{a_1, \cdots, a_i \geq 1}{\sum_{j=1}^i a_i =r}} \prod_{j=1}^i 
\frac{1}{a_j}\cdot \left( \left(1+\frac{2}{A}\right)\cdot q^{a_j} +1\right)
  \\
& \leq 
\sum_{i=1}^r \frac{1}{i!}\sum_{\genfrac{.}{.}{0pt}{}{a_1, \cdots, a_i \geq 1}{\sum_{j=1}^i a_i =r}}
\prod_{j=1}^i 2a \cdot \left(1+ \frac{2}{A} \right) \cdot q^{a_j}
\notag \\
& = q^r\cdot   \sum_{i=1}^r \sum_{\genfrac{.}{.}{0pt}{}{a_1, \cdots, a_i \geq 1}{\sum_{j=1}^i a_i =r}} \frac{1}{i!}\cdot  \left(2a\right)^i \cdot \left(1+ \frac{2}{A} \right)^i \notag  \\
& =
q^r \cdot  \left(2a+ \frac{4a}{A} \right) \cdot  \sum_{i=1}^r  \frac{2^{i-1}}{i!} \cdot \binom{r-1}{i-1}\cdot\left(a+ \frac{2a}{A} \right)^{i-1}.
\notag \\
& \leq 
q^r \cdot  \left(2a+ \frac{4a}{A} \right) \cdot \sum_{i=1}^r \binom{r-1}{i-1} \cdot\left(a+ \frac{2a}{A} \right)^{i-1} \notag \\  
&= q^r \cdot  \left(2a+ \frac{4a}{A} \right) \cdot \left(a+1+ \frac{2a}{A} \right)^{r-1} \notag \\
& < q^r \cdot \left(\frac{5}{3}+ \frac{4}{3A} \right)^{r}. \notag
\end{align}
This completes the proof of the lemma.
\end{proof}
\vspace{3mm}

Next, denote by $J$ the Jacobian variety of $X$,
and 
let us fix 
   a $k$-rational point $P$ of $X$.
    (If the condition $(**)$ is satisfied, then it suffices to choose $P$ from  $R\cup S\cup T$ . In fact, the assumptions $n \geq 1$ and $n \geq g + \mr{max}\{t, g\}$ implies  $\sharp (R \cup S\cup T) = n+s-g \geq n-g \geq 1$.)
In what follows, if $\mcL$ is  a line bundle of degree $0$ on $X$, we abuse notation and also write $\mcL$ for the corresponding point of $J$.
For each positive integer $r$, we denote by
\begin{align}
\varphi_r : X^{(r)} \migi J
\end{align}
the morphism giving by sending each divisor $D$
to the line bundle $\mcO (-rP+ D)$.
It follows from  ~\cite{Mil}, \S\,5, Theorem 5.1 (a), that $\varphi_m$ is surjective (resp., birational onto its image; resp., a birational morphism) if $m > g$ (resp., $m < g$; resp., $m=g$).
The fiber $\varphi^{-1}_r$ of 
this morphism 
over $\mcL \in J (k)$ can be identified with the complete linear system $|\mcL(rP)|$ of $\mcL (rP)$.
In particular, if $k=\mbF_q$ (i.e., the assumption $(**)$ holds), then
the number of $\mbF_q$-rational points of $\varphi^{-1}_r (\mcL)$ is given by
\begin{align} \label{E034}
\sharp \varphi_r^{-1}(\mcL)(\mbF_q) \ \left(= \sharp |\mcL (r P)| (\mbF_q)\right) = 
\sum_{j=0}^{L}q^{j}
 = q^{L} + \frac{q^{L}-1}{q-1} < q^{L}+ \frac{q^{L}}{2} = \frac{3}{2} \cdot q^{L},
\end{align}
where $L := \mr{dim}|\mcL (rP)|$.

Now, we shall 
 denote by
\begin{align}
\iota : X^{(g)} \migi X^{(n)}
\end{align}
the closed immersion given by sending each divisor $D$  to $D +  \sum_{Q \in R \cup T} Q$.
The composite
$\varphi_n \circ \iota : X^{(g)} \migi J$ coincides with the composite of  $\varphi_g$ and  the translation of $J$ by the line bundle $\mcO (-(n-g)P +\sum_{Q \in R \cup T}  Q)$.
It follows that $\varphi_n \circ \iota$ is a  birational morphism.
More precisely, (according to ~\cite{Mil}, Lemma 5.2 (b), and the comment following that lemma) we can find a dense open subset $U$ of $X^{(g)}$ such that
$h^0 (\mcO (D)) =1$ (or equivalently, $h^0 (\Omega_{X/k}(-D))=0$ by the Riemann-Roch theorem) for all $D$'s in $U$; the restriction of $\varphi_n \circ \iota$ to $U$ is an open immersion.
We shall write
\begin{align}
J' := (\varphi_n \circ \iota) (U), \hspace{5mm}
\overline{E}_1 := J \setminus J', \ \ \   \text{and} \ \ \
 E_1 := \varphi^{-1}_n(\overline{E}_1),
 \end{align}
where  $\overline{E}_1$ and $E_1$ are considered as  reduced closed subschemes of $J$ and $X^{(n)}$ respectively.
Since the restriction $\varphi_n |_{\iota (U)} : \iota (U) \migi J'$ of $\varphi_n$ to $\iota (U)$ is an isomorphism, 
there exists its inverse morphism
\begin{align}
\psi : J' \isom \iota (U) \ \left( \subseteq  X^{(n)}\right).
\end{align}

Also, for each $k$-rational point $Q$ of $X$, denote by
\begin{align}
\gamma_Q : X^{(g-1)} \migi X^{(g)} 
\end{align}
the morphism given by $D \mapsto D + Q$.
Then, we have a reduced closed subscheme
\begin{align}
E_{Q} := \varphi^{-1}_n(\mr{Im}(\varphi_n \circ \iota \circ \gamma_Q))
\end{align}
of $X^{(n)}$.

\vspace{3mm}
\ble \label{L001} \leavevmode\\
 \ \ \ 
Let $Q$ be an $k$-rational point of $X$.
Then, 
both $E_1$ and $E_Q$ are of dimension $\leq n-1$.
If, moreover, the condition $(**)$ is satisfied, then
the following inequalities hold:
\begin{align}
\sharp E_1 (\mbF_{q}), \sharp E_Q (\mbF_q) \leq \frac{3}{2} \cdot \left(\frac{5}{3}+ \frac{4}{3A}\right)^{g-1} \cdot  q^{n-1}.
\end{align}
\ele
\begin{proof}
We will only consider the case of $E_1$
    because the proof of the other is similar.
Denote by $\Theta'$ the reduced closed subscheme of $X^{(g)}$ classifying divisors $D$ with $h^0 (\Omega_{X/k}(-D+P))>0$.
Since $\mr{deg}(\Omega_{X/k}(-D+P))=g-1$ for any such  divisor $D$,  it follows from a well-known fact that the image $(\varphi_n \circ \iota) (\Theta')$ coincides with 
 the theta divisor $\Theta := \mr{Im}(\varphi_{g-1}) \ (\subseteq J)$  up to translation.
The obvious inequality $h^0 (\Omega_{X/k}(-D+P)) \geq h^0 (\Omega_{X/k}(-D))$
implies that $(\varphi_n \circ \iota) (\Theta')$ contains $\overline{E}_1$, and hence,
\begin{align} \label{E778}
\mr{dim}(\overline{E}_1) \leq \mr{dim}((\varphi_n \circ \iota) (\Theta')) \leq \mr{dim}(\Theta') = g-1.
\end{align}
Let us take a line bundle $\mcL$ classified by  $\overline{E}_1  \ \left(\subseteq J \right)$.
Since $\mr{deg} (\Omega_{X/{k}}\otimes \mcL (nP)^\vee) = 2g-2-n <0$ ($\Longrightarrow h^1 (\mcL (nP))=0$) by assumption, the Riemann-Roch theorem gives 
\begin{align}
h^0 (\mcL(nP)) = h^0 (\mcL (nP)) - h^1(\mcL (nP)) = n +1-g.
\end{align}
Thus, the equality $\mr{dim}(|\mcL nP|)= n-g$ holds, and 
\begin{align} \label{E777}
\mr{dim}(E_1)  \leq \mr{dim}(|\mcL nP|) + \mr{dim}(\overline{E}_1) \stackrel{(\ref{E778})}{\leq} (n-g) + (g-1) = n-1
\end{align}
(cf. ~\cite{Har}, Chap.\,II, \S\,3, Exercise 3.22 (b)).  

Next,  we shall estimate the number $\sharp E_1 (\mbF_q)$ when the condition $(**)$ holds.
By the comment preceding (\ref{E778}), we have
\begin{align} \label{E017}
\sharp \overline{E}_1 (\mbF_q) \leq \sharp (\varphi_n \circ \iota) (\Theta')(\mbF_q) = \sharp \Theta (\mbF_q).
\end{align}
Moreover, the fiber of $\varphi_{g-1}  : X^{(g-1)} \migi (\Theta \subseteq) \  J$ 
over each $\mbF_q$-rational point of $\Theta$
is isomorphic to a projective space (which, in particular,  contains at least one $\mbF_q$-rational point).
This implies  that $\varphi_g$ induces a surjective map
$X^{(g-1)} (\mbF_q) \migisurj \Theta (\mbF_q)$, and hence,
\begin{align} \label{E020}
\sharp \Theta (\mbF_q) \leq  \sharp X^{(g-1)} (\mbF_q).
\end{align}
By (\ref{E017}), (\ref{E020}), and Lemma \ref{L005} (ii), the following sequence of inequalities holds:
\begin{align} \label{E015}
\sharp \overline{E}_1 (\mbF_q) \leq \sharp \Theta (\mbF_q)\leq \sharp X^{(g-1)} (\mbF_q) <   \left(\frac{5}{3}+ \frac{4}{3A}\right)^{g-1} \cdot  q^{g-1}.
\end{align}
On the other hand, for each line bundle $\mcL$ classified by $\overline{E}_1$, the  result of (\ref{E034}) reads
\begin{align}\label{E014}
\sharp \varphi^{-1}_n (\mcL) (\mbF_q) < \frac{3}{2} \cdot q^{n-g}.
\end{align}
Thus, by (\ref{E015}) and (\ref{E014}), we have
\begin{align}
\sharp E_1 (\mbF_{q}) 
&\leq \sum_{\mcL \in \overline{E}_1 (\mbF_q)}\sharp \varphi_n^{-1}(\mcL)(\mbF_q)
\\
&\leq   \left(\frac{5}{3}+ \frac{4}{3A}\right)^{g-1} \cdot q^{g-1}  \cdot \frac{3}{2}\cdot q^{n-g} \notag \\
&= \frac{3}{2} \cdot \left(\frac{5}{3}+ \frac{4}{3A}\right)^{g-1} \cdot  q^{n-1}.\notag
\end{align}
This completes the proof of the lemma.
\end{proof}
\vspace{3mm}

Next, let us consider 
the reduced closed subscheme 
\begin{align}
E_2 := \varphi^{-1}_n (\varphi_n (\mr{Im}(\alpha^{(n)}) \cap \mr{Im}(\iota)))
\end{align}
of $X^{(n)}$, where, for each positive integer $r$, $\alpha^{(r)}$ denotes the morphism $X^{r-1} \migi X^{(r)}$ given by
\begin{align}
\alpha^{(r)} (Q_1, \cdots, Q_{r-1}) &:= 2Q_1 + \sum_{i=2}^{r-1}Q_i.
\end{align}

\vspace{3mm}
\ble \label{L0012} \leavevmode\\
 \ \ \ 
$E_2$ is of dimension $\leq n-1$.
If, moreover,
the condition $(**)$ is satisfied,
  then
 the following inequality holds:
 \begin{align}
 \sharp E_2 (\mbF_q) < \left(\frac{3n}{2}+\frac{3g}{8}\right)  \cdot \left(\frac{5}{3}+ \frac{4}{3A}\right)^{g-1} \cdot q^{n-1}.
 \end{align}
\ele
\begin{proof}
Each divisor classified by a point of $\mr{Im}(\alpha^{(n)}) \cap \mr{Im}(\iota))$ may be expressed as  $D +  \sum_{Q \in R \cup T} Q$, where $D$ is a divisor in either $\mr{Im}(\iota \circ \alpha^{(g)})$ or $\mr{Im}(\iota \circ \gamma_Q)$ for some $Q \in R \cup T$.
This implies that
\begin{align}
\mr{dim}(\varphi_n(\mr{Im}(\alpha^{(n)}) \cap \mr{Im}(\iota))))  
&\leq \mr{dim}(\mr{Im}(\alpha^{(n)}) \cap \mr{Im}(\iota)))   \\
 &\leq \mr{max} \left\{\mr{dim}(\mr{Im}(\iota \circ \alpha^{(g)})), \mr{dim}(\mr{Im}(\iota \circ \gamma_Q)) \right\} \notag \\
 & \leq  \mr{max} \left\{ \mr{dim}(\mr{Im}(\alpha^{(g)}))\right\} \cup \left\{ \mr{dim}(\mr{Im}(\gamma_Q))\right\}_{Q \in R \cup T} \notag \\
 & \leq \mr{max} \left\{ \mr{dim} (X^{g-1}), \mr{dim}(X^{(g-1)}) \right\} \notag \\
 & = g-1. \notag
\end{align}
Hence, by the same argument as the inequalities (\ref{E777}), we have
\begin{align}
\mr{dim}(E_2) \leq  (n-g)+\mr{dim}(\varphi_n(\mr{Im}(\alpha^{(n)}) \cap \mr{Im}(\iota)))) = n-1,
\end{align}
which completes the proof of the former assertion.

Next, let us consider the latter assertion.
Suppose that the condition $(**)$ is satisfied.
Then, the above discussion implies that
\begin{align} \label{E035}
\sharp \mr{Im}(\alpha^{(n)}) \cap \mr{Im}(\iota)) (\mbF_q)
&\leq
\sharp  \mr{Im}(\iota \circ \alpha^{(g)} ) (\mbF_q)+  \sum_{Q \in R \cup T} \sharp \mr{Im}(\iota \circ \gamma_Q)(\mbF_q)  \\
& \leq \sharp \mr{Im}(\alpha^{(g)}) + \sum_{Q \in R \cup T} \sharp X^{(g-1)}(\mbF_q) \notag \\
& \leq \sharp \mr{Im}(\alpha^{(g)}) + (n-g) \cdot \left(\frac{5}{3}+ \frac{4}{3A} \right)^{g-1} \cdot q^{g-1},\notag
\end{align}
where the last inequality follows from Lemma \ref{L005} (ii).

In what follows, let us consider an upper bound for $\sharp \mr{Im}(\alpha^{(g)})$.
For each  $\overline{\mbF}_q$-rational point $Q$ of $X$,
we shall denote by $Q_{\mr{div}}$ the reduced effective divisor determined by its image.
The assignment $Q \mapsto Q_\mr{div}$ can be extended naturally to a map $D \mapsto D_\mr{div}$ from   $X^{(r)}(\overline{\mbF}_q)$ (for each $r$) to the set of divisors on $X$.
Now, let us write
\begin{align} \label{E040}
 \breve{X}^{(m; l)} := \left\{ (D, D') \in X^{(m)} (\mbF_{q}) \times X^{(l)}(\overline{\mbF}_q) \, | \,  D \geq  D'_\mr{div} \right\}
\end{align}
($m \geq l \geq 1$).
Since  the fiber of the first projection  $\breve{X}^{(m; l)} \migi  X^{(m)} (\mbF_{q})$ over each element of $X^{(m)} (\mbF_{q})$ has at most $\binom{m}{l}$ elements.
It follows that
\begin{align} \label{E038}
\breve{X}^{(m; l)} \leq \binom{m}{l}  \cdot X^{(m)} (\mbF_{q}) \leq  \binom{m}{l} \cdot \left(\frac{5}{3}+\frac{4}{3A} \right)^m \cdot q^m.
\end{align}
Now, let us consider the set
\begin{align}
\breve{X}_2^{(g)} := \left\{ (D, D') \in \mr{Im}(\alpha^{(g)}) (\mbF_{q}) \times X(\overline{\mbF}_q) \, | \,  D \geq 2 Q_\mr{div} \right\}.
\end{align}
The projection to the first factor
$\breve{X}_2^{(g)} \migi \mr{Im}(\alpha^{(g)}) (\mbF_{q})$
is surjective, and hence we have
\begin{align} \label{E043}
\sharp  \mr{Im}(\alpha^{(g)}) (\mbF_{q}) \leq  \sharp \breve{X}_2^{(g)}.
\end{align}

Moreover, notice that the assignment $(D, D') \mapsto D-D'_\mr{div}$ give
an injection
\begin{align}
\breve{X}_2^{(g)} \migiincl \coprod_{i=1}^{g-1} \breve{X}^{(i; 1)}.
\end{align}
Hence, by (\ref{E038}), we have
\begin{align} \label{E045}
 \breve{X}_2^{(g)} 
&\leq \sum_{i=1}^{g-1} \breve{X}^{(i; 1)}  \\
&\leq \sum_{i=1}^{g-1} i  \cdot \left(\frac{5}{3}+ \frac{4}{3A} \right)^i \cdot q^i \notag \\
&< g \cdot  \left(\frac{5}{3}+\frac{4}{3A} \right)^{g-1}\cdot q^{g-1}  + g  \cdot  \left(\frac{5}{3}+ \frac{4}{3A} \right) \cdot q \cdot  \frac{( (5/3+4/3A) \cdot q)^{g-2}-1}{ (5/3+4/3A) \cdot q-1} \notag \\
& < g \cdot \left(\frac{5}{3}+\frac{4}{3A} \right)^{g-1}\cdot q^{g-1}  + g   \cdot \left(\frac{5}{3}+ \frac{4}{3A} \right) \cdot q \cdot  \frac{1}{4} \cdot  \left(\frac{5}{3}+ \frac{4}{3A} \right)^{g-2} \cdot q^{g-2}\notag \\
& = \frac{g \cdot 5}{4}  \cdot \left(\frac{5}{3}+ \frac{4}{3A} \right)^{g-1} \cdot q^{g-1}. \notag
\end{align} 

By combining (\ref{E035}), (\ref{E043}), and (\ref{E045}), we obtain
\begin{align} \label{E050}
\sharp \mr{Im}(\alpha^{(g)}\cap \mr{Im}(\iota))(\mbF_q) 
&<
\frac{g \cdot 5}{4}  \cdot \left(\frac{5}{3}+ \frac{4}{3A} \right)^{g-1} \cdot q^{g-1} + (n-g) \cdot \left(\frac{5}{3}+ \frac{4}{3A} \right)^{g-1} \cdot q^{g-1}
  \\
&= \left(n+\frac{g}{4}\right)  \cdot \left(\frac{5}{3}+ \frac{4}{3A} \right)^{g-1} \cdot q^{g-1}. \notag
\end{align}
Finally, since each fiber of $\varphi_n$ is isomorphic to a projective space, 
the following inequalities hold:
\begin{align}
\sharp E_2 (\mbF_q)  
&\leq \sum_{\mcL \in \varphi_n(\mr{Im}(\alpha^{(n)}\cap \mr{Im}(\iota))(\mbF_q))} \sharp \varphi^{-1}_n (\mcL)(\mbF_q)  \\
& \leq 
\sum_{D \in \mr{Im}(\alpha^{(n)}\cap \mr{Im}(\iota))(\mbF_q)} \sharp \varphi^{-1}_n (\varphi_n (D))(\mbF_q) \notag \\
& \leq \sum_{D \in \mr{Im}(\alpha^{(n)}\cap \mr{Im}(\iota))(\mbF_q)}
\frac{3}{2} \cdot q^{n-q} \notag \\
& <  \left(n+\frac{g}{4}\right)  \cdot \left(\frac{5}{3}+ \frac{4}{3A} \right)^{g-1}  \cdot q^{g-1} \cdot \frac{3}{2} \cdot q^{n-q} 
 \notag \\
&< \left(\frac{3n}{2}+\frac{3g}{8}\right)  \cdot \left(\frac{5}{3}+ \frac{4}{3A}\right)^{g-1} \cdot q^{n-1}, \notag 
\end{align}
where the third inequality  follows from (\ref{E034}) and the last inequality follows from (\ref{E050}).
This completes the proof of the lemma.
\end{proof}
\vspace{3mm}

Let 
$\beta : X^{n-2} \migi X^{(n)}$ and $\gamma : X^{n-2} \migi X^{(n)}$ be  the morphisms given by 
\begin{align}
\beta (Q_1, \cdots, Q_{n-2}) := 2Q_1+2Q_2+ \sum_{l=3}^{n-2} Q_l, \hspace{5mm}
\gamma (Q_1, \cdots,  Q_{n-2}) := 3Q_1 + \sum_{l=2}^{n-2} Q_l.
\end{align}
Write
\begin{align}
\delta : ((\mr{Im}(\beta) \cup \mr{Im}(\gamma)) \setminus E_1) \times \mbP^1_k \migi X^{(n)}
\end{align}
given by $(D, t) \mapsto \psi \circ \varphi_n (D) + t (D- \psi \circ \varphi_n (D))$, and write
\begin{align}
E_3 := \mr{Im}(\delta)
\end{align}
(considered as a reduced  subscheme of $X^{(n)}$).

\vspace{3mm}
\ble \label{L002} \leavevmode\\
 \ \ \ 
$E_3$ is of dimension $\leq n-1$.
If, moreover,
the condition $(**)$ is satisfied,
 then
the following inequality holds:
 \begin{align}
 \sharp E_3 (\mbF_q) < \frac{5n^2}{2}  \cdot \left(\frac{5}{3}+ \frac{4}{3A} \right)^{n-2}  \cdot q^{n-1}.
 \end{align}
\ele
\begin{proof}
Since $\mr{dim}(\mr{Im}(\beta))$, $\mr{dim}(\mr{Im}(\gamma)) \leq \mr{dim}(X^{n-2}) = n-2$,
the following inequalities hold:
\begin{align}
\mr{dim}(E_3)\leq \mr{dim}(\mbP^1_k) + \mr{max}\left\{\mr{dim}(\mr{Im}(\beta)),\mr{dim}(\mr{Im}(\gamma)) \right\}\leq 1+(n-2)=n-1.
\end{align}
Next, we shall prove the latter assertion under the assumption that $(**)$ is satisfied.
Let us write
\begin{align}
\breve{X}_\beta^{(n)} &:= \left\{ (D, D') \in \mr{Im}(\beta) (\mbF_{q}) \times X^{(2)}(\overline{\mbF}_q) \, | \,  D \geq 2 Q_\mr{div} \right\},  \\
\breve{X}_\gamma^{(n)} &:= \left\{ (D, D') \in \mr{Im}(\gamma) (\mbF_{q}) \times X(\overline{\mbF}_q) \, | \,  D \geq 3 Q_\mr{div} \right\}.\notag
\end{align}
The projections to the first factors  $\breve{X}_\beta^{(n)} \migi \mr{Im}(\beta)(\mbF_q)$ and $\breve{X}_\gamma^{(n)} \migi \mr{Im}(\gamma)(\mbF_q)$ are surjective, and we have
\begin{align} \label{E051}
\sharp  \mr{Im}(\beta)(\mbF_q) \leq \sharp \breve{X}_\beta^{(n)}, \hspace{5mm}  \sharp \mr{Im}(\gamma)(\mbF_q) \leq  \breve{X}_\gamma^{(n)}.
\end{align}
Moreover, notice that the assignments  $(D, D') \mapsto D-D'_\mr{div}$ and $(D, D') \mapsto D-2D'_\mr{div}$ give
injections 
\begin{align} \label{E052}
\breve{X}_\beta^{(n)} \migiincl \coprod_{i=2}^{n-2} \breve{X}^{(i; 2)},
 \hspace{5mm}
 \breve{X}_\gamma^{(n)} \migiincl \coprod_{i=1}^{n-2} \breve{X}^{(i; 1)}
\end{align}
respectively, where $\breve{X}^{(m; l)}$'s are as in the proof of Lemma \ref{L0012}.
Hence, by (\ref{E051}), (\ref{E052}), and (\ref{E038}), the following inequalities hold:
\begin{align} \label{E025}
\sharp \mr{Im}(\beta) &\leq \sum_{i=2}^{n-2} \breve{X}^{(i; 2)} \\
& \leq \sum_{i=2}^{n-2} \binom{i}{2} \cdot \left(\frac{5}{3}+ \frac{4}{3A} \right)^{i} \cdot q^i \notag  \\
& \leq \frac{n^2}{2} \cdot \left(\frac{5}{3}+ \frac{4}{3A} \right)^{n-2} \cdot q^{n-2} + \frac{n^2}{2} \cdot \left(\frac{5}{3}+ \frac{4}{3A} \right)^{2} \cdot q^2\cdot \frac{\left((5/3+4/3A)q\right)^{n-4}-1}{(5/3+4/3A)q-1} \notag  \\
& < \frac{n^2}{2} \cdot \left(\frac{5}{3}+ \frac{4}{3A} \right)^{n-2} \cdot q^{n-2} + \frac{n^2}{2} \cdot \frac{1}{4} \cdot \left(\frac{5}{3}+ \frac{4}{3A} \right)^{n-2} \cdot q^{n-2}  \notag \\
& = \frac{5n^2}{8} \cdot \left(\frac{5}{3}+ \frac{4}{3A} \right)^{n-2} \cdot q^{n-2}, \notag
 \\
\sharp \mr{Im}(\gamma) &\leq \sum_{i=1}^{n-2} \breve{X}^{(i; 1)} \\
& \leq \sum_{i=1}^{n-2} i \cdot \left(\frac{5}{3}+ \frac{4}{3A} \right)^{i} \cdot q^i \notag \\
& <  n  \cdot \left(\frac{5}{3}+ \frac{4}{3A} \right)^{n-2} \cdot q^{n-2} + n \cdot \left(\frac{5}{3}+ \frac{4}{3A} \right) \cdot q \cdot \frac{\left((5/3+4/3A)q\right)^{n-3}-1}{(5/3+4/3A)q-1} \notag \\
& < n  \cdot \left(\frac{5}{3}+ \frac{4}{3A} \right)^{n-2} \cdot q^{n-2} + 
n\cdot \frac{1}{4} \cdot \left(\frac{5}{3}+ \frac{4}{3A} \right)^{n-2} \cdot q^{n-2}
\notag \\ 
& = \frac{5n}{4} \cdot \left(\frac{5}{3}+ \frac{4}{3A} \right)^{n-2}  \cdot q^{n-2}. \notag
\end{align} 
Thus, we have
\begin{align}
\sharp E_3 (\mbF_q) 
& \leq 
(\sharp \mr{Im}(\beta) + \sharp \mr{Im}(\gamma)) \cdot \sharp \mbP^1 (\mbF_q) \\
& < \left(\frac{5n^2}{8} \cdot \left(\frac{5}{3}+ \frac{4}{3A} \right)^{n-2}  \cdot q^{n-2}+ \frac{5n}{4} \cdot \left(\frac{5}{3}+ \frac{4}{3A} \right)^{n-2}  \cdot q^{n-2}\right) \cdot (q+1) \notag \\
& = \frac{5n (n+2)}{8}  \cdot \left(\frac{5}{3}+ \frac{4}{3A} \right)^{n-2}  \cdot q^{n-2} \cdot (q+1) \notag \\
&< \frac{5n^2}{2}  \cdot \left(\frac{5}{3}+ \frac{4}{3A} \right)^{n-2}  \cdot q^{n-1}.  \notag 
\end{align}
This completes the proof of the lemma.
\end{proof}
\vspace{3mm}

Denote by
\begin{align}
E'_4
\end{align}
the closed subscheme of $X^{g} \times_{\varphi_n \circ \iota \circ \pi_g, J, \varphi_n} X^{(n)}$ consisting of points $((Q_i)_{i=1}^g), D)$ 
such that $Q \in \mr{Supp} (D)$ for some $Q \in R \cup T$ or $Q = Q_i$ ($i=1, \cdots, g$).
Denote by 
\begin{align}
E_4
\end{align}
the scheme-theoretic image of the second projection $E'_4 \migi X^{(n)}$.

\vspace{3mm}
\ble \label{L008} \leavevmode\\
 \ \ \ 
$E_4$ is of dimension $\leq n-1$.
If, moreover,
the condition $(**)$ is satisfied, then the following inequality holds:
\begin{align}
\sharp E_4 (\mbF_q) <\frac{3}{2} \cdot \left(\frac{5}{3} + \frac{4}{3A} \right)^g \cdot q^{n-1}.
\end{align}
\ele
\begin{proof}
Let $D \in X^{(n)}$.
Since $\mr{deg}(D) = n \geq 2g$, the linear system $|D|$ has no base points (cf. ~\cite{Har}, Chap.\,IV, \S\,3, Corollary 3.2).
This implies that $|D(-Q)| \neq |D|$ (i.e., $\mr{dim}(|D (-Q)|) \leq n-g-1$) for  each point $Q$ of $X$.
Hence, 
 each fiber of the morphism  $\varphi_n |_{E_4} : E_4 \migi J$ induced by $\varphi_n$ is contained in the union of $n$ hyperplanes embedded in the $n-g$ dimensional projective space. 
By the same argument as the inequalities (\ref{E777}),
we have
\begin{align}
\mr{dim}(E_4) \leq \mr{dim}(E'_4) \leq \mr{dim}(J) + (n-g-1) = n-1,
\end{align}
which completes the former assertion.

Also, if the condition $(**)$ is satisfied,
then the above discussion shows that
$\sharp (\varphi_n |_{E_4})^{-1}(\mcL) \leq n \cdot \frac{3}{2} \cdot q^{n-g-1}$ for each $\mcL \in J (\mbF_q)$ (cf.  (\ref{E034})).
  Thus, 
\begin{align}
\sharp E_4 (\mbF_q) &\leq 
\sum_{\mcL \in J (\mbF_q)}  \sharp (\varphi_n |_{E_4})^{-1}(\mcL) \\
&\leq \sharp X^{(g)}(\mbF_q)  \cdot \left(n \cdot \frac{3}{2} \cdot q^{n-g-1}\right) \notag\\
& < \left(\frac{5}{3} + \frac{4}{3A} \right)^g \cdot q^{g} \cdot \left(n \cdot \frac{3}{2} \cdot q^{n-g-1}\right)\notag \\
& \leq \frac{3}{2} \cdot \left(\frac{5}{3} + \frac{4}{3A} \right)^g \cdot q^{n-1}. \notag
\end{align}
This completes the proof of the lemma.
\end{proof}
\vspace{3mm}

Now, let us complete the proof of Proposition \ref{P02}.
Write $E := E_1 \cup E_2 \cup E_3 \cup E_4 \cup \bigcup_{Q \in S}E_Q$.
If the condition $(**)$ is satisfied,  then the lemmas proved so far and the assumption (\ref{Erfhih}) imply the following sequence of inequalities: 
 \begin{align}
 & \ \ \ \ \sharp (X^{(n)}\setminus E)(\mbF_q) \\
 & \geq \sharp X^{(n)}(\mbF_q)  - \sharp E_1 (\mbF_q) - \sharp E_2 (\mbF_q) - \sharp E_{3} (\mbF_q) - \sharp E_{4} (\mbF_q) - \sum_{Q \in S} \sharp E_{Q} (\mbF_q) 
\notag \\
 & > \frac{1}{7 \cdot n!} \left(3-\frac{2}{A} \right)^n \cdot q^n  
 - \frac{3}{2} \cdot \left(\frac{5}{3} +\frac{4}{3A}\right)^{g-1} \cdot q^{n-1}
 - \left(\frac{3n}{2}+ \frac{3g}{8} \right) \cdot \left(\frac{5}{3} + \frac{4}{3A} \right)^{g-1}  \cdot q^{n-1} \notag  \\
 & \ \ \ \ -\frac{5n^2}{2}  \cdot \left(\frac{5}{3}+ \frac{4}{3A} \right)^{n-2}\cdot q^{n-1} - \frac{3}{2} \cdot \left(\frac{5}{3}+ \frac{4}{3A} \right)^{g} \cdot q^{n-1}- s \cdot \frac{3}{2} \cdot  \left(\frac{5}{3}+ \frac{4}{3A} \right)^{g-1} \cdot q^{n-1} \notag \\
  & \geq \frac{1}{7 \cdot n!} \left(\frac{3A-2}{A} \right)^n \cdot q^{n-1} \left(q- 10^2 \cdot n! \cdot  (n^2 +s) \cdot \left(\frac{5A+4}{9A-6} \right)^{n} \right) \notag \\
  & > 0. \notag
 \end{align}
Hence, 
 there exists a $k$-rational point 
 $Q$ of $X^{(n)}$ which is not contained in $E$.
We can also verify this claim in the case of $(*)$ because of  the dimension estimates proved 
 in the lemmas, which implies that $\mr{dim}(E)\leq n-1$.
 Denote by $l$ the line in the linear system  $|\mcO_X (D) |$ passing through $Q$ and $(\psi \circ \varphi_n)(Q)$. (Notice that $(\psi \circ \varphi_n)(Q)$ can be defined since $Q$ is not contained in $E_1$.)
  Choose $\zeta_0$, $\zeta_1 \in H^0 (X, \mcO_X (D))$ such  that $l$ consists of all divisors $\mr{div} (\lambda_1 \zeta_0 - \lambda_0 \zeta_1)$, where $\lambda_0$, $\lambda_1 \in k$ are not both zero.
  Since $D$ is not contained in  $E_4$,
  two divisors corresponding to $Q$ and $(\psi \circ \varphi_n)(Q)$ have no common support.
Hence,  the morphism $\zeta : X \migi \mbP^1_k$
 given by $\zeta (P) = [\zeta_0 (P): \zeta_1 (P)]$ is everywhere-defined and of degree $n$.
In particular, $\zeta^{-1}((\lambda_0, \lambda_1)) = \mr{div} (\lambda_1 \zeta_0 - \lambda_0 \zeta_1)$.
Since $l \cap \bigcup_{Q \in S}E_Q =\emptyset$ (resp., $l \cap E_2 =\emptyset$; resp., $l \cap E_3 = \emptyset$), we see that $\zeta (T) \cap \zeta (S) = \emptyset$ (resp., $\zeta$ is unramified over $\zeta (T)$; resp., $\zeta$ is a simple covering by ~\cite{Ful}, Theorem 5.6).
Consequently, $\zeta$ specifies the desired morphism.
This completes the proof of Proposition \ref{P02}. 
 \end{proof}
\vspace{3mm}

\vspace{10mm}
\section{Second step of the proof} \vspace{3mm}

Next, we prove the following effective version of ~\cite{SZ}, Proposition 4.1.

\vspace{3mm}
\bpr \label{L04} \leavevmode\\
 \ \ \ 
Let $q$ be a power of an odd prime, and let $S$ be a (possibly empty) finite set of $\mbF_q$-rational points of $\mbP^1_{\mbF_q}$.
Write $s := \sharp S$.
\begin{itemize}
\item[(i)]
There exists a $p$-tame Bely\I \ map on $(\mbP^1_{\mbF_q}, \emptyset, S)$ over $\mbF_q$ of degree $q-1$.
\item[(ii)]
Let $\tau$ be an $\mbF_q$-rational point of $\mbP^1_{\mbF_q}$.
If $s \leq 3$ (resp., $s > 3$), then there exists a $p$-tame Bely\I \ map on $(X, S, \{ \tau \})$ over $\mbF_q$ of degree $1$ (resp., $(q-1)^{s-3}$).
\end{itemize}
 \epr
\begin{proof}
Assertion (i) follows from the fact that
the $\mbF_q$-endomophism of the field $\mbF_q (x)$ (= the rational function field of one variable $x$ over $\mbF_q$)  given by $x \mapsto x^{q-1}-1$ specifies the desired endomorphism of $\mbP^1_{\mbF_q}$ (since any element $v$ of $\mbF_q  = \mbP^1_{\mbF_q} (\mbF_q) \setminus \{ \infty \}$ satisfies the equality $v^{q-1} =1$).

Next, we shall consider assertion (ii).
The non-resp'd portion is immediately verified by taking a suitable linear transformation on $\mbP^1_{\mbF_q}$.
In what follows, let us prove the resp'd portion 
 by induction on $s$.
After possibly applying some linear transformation on $\mbP^1_{\mbF_q}$, we can suppose  that 
 $\{0, \infty \} \subseteq  S$.
Let us fix $\alpha \in S \setminus \{0, \infty \}$, and
denote by
$\xi_1$ 
the $\mbF_q$-endomorphism of $\mbP^1_{\mbF_q}$ (of degree $q-1$) corresponding to the $\mbF_q$-endomorphism of $\mbF_q (x)$
  given by $x \mapsto -x^{q-1}+ \alpha^{-1}x$.
  Then, $\xi_1 (0) = \xi_1 (\alpha) = 0$, and $\xi_1$ is unramified away from $\infty$ (since $(2, q) = 1$).
 Since  $\xi_1 (\beta) = \alpha^{-1}\beta -1$ for any $\beta \in \mbF_q$,
  the map of sets $\mbF_p \setminus \{ 0, \alpha\} \migi \mbF_q \setminus \{0\}$ induced by
  $\xi_1$ is injective.
Hence, $\sharp \xi_1 (S) = s-1$ and $\xi_1 (\tau) \notin \xi_1 (S)$.
By applying the inductive assumption (or the non-resp'd portion already proved),
we obtain a $p$-tame Bely\I \ map $\xi_2 : \mbP^1_{\mbF_q} \migi \mbP^1_{\mbF_q}$ on $(\mbP^1_{\mbF_q}, \psi_1 (S), \{ \xi_1 (\tau) \})$ over $\mbF_q$ of degree $(q-1)^{s -4}$.
The composite $\xi := \xi_2 \circ \xi_1$ has degree $\mr{deg}(\xi_2) \cdot \mr{deg}(\xi_1) = (q-1)^{s -4} \cdot (q-1) = (q-1)^{s -3}$ and gives the desired morphism.
This completes the proof of the proposition.
\end{proof}
\vspace{3mm}

\begin{proof}[Proof of Theorem A]
By the last assertion of Proposition \ref{P02},
we see that if $m$ is a positive integer with $q^m \geq 10^2 \cdot (2g+t+1)! \cdot (2g+t+s+1)^2  \cdot (5/6)^{2g+t+1}$, then there exists a simple covering $\zeta : X_{\mbF_{q^m}} \migi \mbP^1_{\mbF_{q^m}}$ over $\mbF_{q^m}$ of degree $2g+t+1$ such that $\zeta (S_{\mbF_{q^m}}) \cap \zeta (T_{\mbF_{q^m}}) = \emptyset$ and $\zeta (T_{\mbF_{q^m}})$ consists of one point $\tau_0$ (if $T \neq \emptyset$).
Indeed, we can choose $m$ as  
\begin{align}
m := \lceil\mr{log}_q (10^2 \cdot (2g+t+1)! \cdot (2g+t+s+1)^2  \cdot (5/6)^{2g+t+1})\rceil.
\end{align}
Since the discriminant $\delta (\zeta)$ of $\zeta$ has degree $2g+2(2g+t+1)-2 = 6g +2t$,
all the points of $\mr{Br} (\zeta)$ are $\mbF_{q^{m \cdot L(6g +2t)}}$-rational.
Hence, if $S'$ denotes the set of closed points of $X_{\mbF_{q^m}}$ defined as the   pull-back of $\zeta (S) \cup \mr{Br}(\zeta)$, then all elements of $S'$ are $\mbF_{q^{m \cdot L(6g +2t)}}$-rational and  
 $\sharp S' \leq  6g+ s+2t$.
Let us choose a $p$-tame Bely\I \ map $\xi$ resulting from   Proposition \ref{L04} (i) or (ii), where the triple $(q, S, \tau)$ is taken to be $(q^{m \cdot L(6g +2t)}, S', \tau_0)$.
Then, the following inequality holds:
\begin{align}
\mr{deg}(\xi) \leq (q^{q^{m \cdot L(6g +2t)}}-1)^{6g+s+2t+1}.
\end{align}
One verifies that the composite $f := \xi \circ \zeta$ 
has the degree
\begin{align}
\mr{deg} (f)  = \mr{deg}(\zeta) \cdot \mr{deg}(\xi)  
 =  (2g+t+1)  \cdot (q^{m \cdot L(6g+2t)}-1)^{6g+s+2t+1},
\end{align}
and moreover, satisfies the required conditions.
This completes the proof of Theorem A. 
\end{proof}

\vspace{10mm}
\section{Appendix: Effective noncritical $p$-wild Bely\I \ theorem} \vspace{3mm}

 The $p$-wild version of our main theorem follows from an argument in the previous works, as we will discuss in the proof of the theorem below.
 Let us keep the notation preceding  Theorem A.
 By a {\bf $p$-wild Bely\I \ map} on $(X, S, T)$ over a field $k$ (where $k$ denotes a field over $\mbF_q$), we shall mean a $k$-dominant morphism $f : X_k \migi \mbP^1_{k}$ satisfying the following conditions:
\begin{align}
f(S_{k}) \cup \mr{Br}(f) \subseteq \{\infty \}, \hspace{10mm} \{ \infty \} \cap f(T_{k}) = \emptyset.
\end{align}
The {\bf $p$-wild Bely\I \ degree} of $(X, S, T)$ is defined as
\begin{align}
{^w \mcB} (X, S, T) := \mr{min} \left\{ \mr{deg}(\phi) \, | \, \text{$\phi$ is  a $p$-wild Bely\I \ map $\phi$ on $(X,S, T)$ over $\overline{\mbF}_q$}\right\}.
\end{align}

\vspace{3mm}
\bt \label{P044} \leavevmode\\
 \ \ \ 
 Suppose that the following inequality holds:
 \begin{align} \label{E090}
q+1 - 2g \sqrt{q}  \geq N +s,
 \end{align}
where  $N := \mr{max}\{2g-1+t, t, 2 \}$.
Then, there exists a $p$-wild Bely\I \ map  on $(X, S, T)$ 
 of degree  
 $< N \cdot p^{s+ 2(g+N)}$.
 In particular, we have
 \begin{align}
 {^w \mcB} (X, S, T) < N \cdot p^{s+ 2(g+N)}.
 \end{align}
 \et
\begin{proof}
By the Hasse-Weil theorem (\ref{E010}) and the 
assumption (\ref{E090}),
  we can find distinct $\mbF_q$-rational points $P_1, \cdots, P_{N-t}$ 
with
 $\{ P_1, \cdots, P_{N-t} \} \cap (S \cup T) = \emptyset$.
Also, it follows from the discussion in ~\cite{SZ}, Lemma 2.2, that there exists  an element 
\begin{align}
v \in \Gamma (X, \mcO_X (\sum_{i=1}^{N-t} P_i + \sum_{Q \in T} Q))\setminus \bigcup_{Q \in T} \Gamma (X, \mcO_X (\sum_{i=1}^{N-t} P_i + \sum_{Q' \in T\setminus \{Q\}} Q')),
\end{align}
and the morphism $\phi : X_{} \migi \mbP^1_{\mbF_q}$  determined by $v^{-1}$ is of degree $N$ and satisfies that 
 $\phi (T) \subseteq \{ 0 \}$ and $\{0 \} \cap (\phi (S)\cup \mr{Br}(\phi))= \emptyset$.
Let
$\overline{B}$ be the set of $\mr{Gal}(\overline{\mbF}_q/\mbF_q)$-conjugates of elements in  $(\phi (S)\cup \mr{Br}(\phi))\setminus \{ \infty \}$, and let
 $V$ denote the $\mbF_p$-span 
  of $\overline{B} \subseteq \mbP^1_{\mbF_q}(\overline{\mbF}_q) \setminus \{\infty \} = \overline{\mbF}_q$.
Also, set
\begin{align}
h_1:= x^p + \frac{x^p}{x^p \cdot h_0(x) + h_0(x)^p}, \hspace{5mm}
h_2 := h_1^p + h_1,
\end{align}
where
$h_0(x) := \prod_{\alpha \in V} (x - \alpha) \in \mbF_q (x)$ (hence $x | h_0 (x)$).
That is to say, $h_1$ and $h_2$ are the rational functions ``$g$" and  ``$f$", respectively,  in the proof of ~\cite{SZ}, Proposition 4.1, where $B$ is taken to be $\phi (S) \cup \mr{Br}(\phi)$.
According to the discussion in {\it loc.\,cit.}, 
the endomorphism $\psi$
  of $\mbP^1_{\mbF_q}$ determined by $h_2$ satisfies that $\mr{Br}(\psi) = \{ \infty \}$ and $\psi (0) \neq \infty$.
Therefore, the composite 
\begin{align}
f := \psi \circ \phi : X \migi \mbP^1_{\mbF_q}
\end{align}
forms a $p$-wild Bely\I \ map on  $(X, S, T)$.
Finally, we shall compute its degree.
By the Riemann-Hurwitz theorem, the discriminant of  $\phi$ turns out to have  degree $2(g +N-1) \ (>1)$, which implies the inequalities $1\leq \sharp (\overline{B}) \leq  s + 2(g+N-1)$.
 Hence, $2<\sharp V \ (= \mr{deg}(h_0)) \leq p^{s + 2(g+N-1)}$,
 and
 \begin{align}
 \mr{deg} (h_2)= p \cdot \mr{deg}(h_1) = p^2 \cdot (\mr{deg}(h_0)-1) < p^{s+ 2(g+N)}.
 \end{align}
 Consequently, 
\begin{align}
\mr{deg}(f) = \mr{deg}(\psi) \cdot \mr{deg}(\phi) <  N \cdot p^{s + 2(g+N)}.
\end{align}
This completes the proof of the proposition.
\end{proof}
\vspace{3mm}

\end{document}